\documentclass[letterpaper]{article} 
\usepackage{aaai2026}  
\usepackage{times}  
\usepackage{helvet}  
\usepackage{courier}  
\usepackage[hyphens]{url}  
\usepackage{graphicx} 
\usepackage{todonotes}
\usepackage{natbib}  
\usepackage{caption} 
\usepackage{url}
\usepackage{algorithm}
\usepackage{algorithmic}
\usepackage{lineno}
\usepackage{newfloat}
\usepackage{listings}
\DeclareCaptionStyle{ruled}{labelfont=normalfont,labelsep=colon,strut=off} 
\floatstyle{ruled}
\newfloat{listing}{tb}{lst}{}
\floatname{listing}{Listing}
\usepackage{complexity}
\usepackage{amsmath}
\usepackage{amssymb}
\usepackage{amsthm}

\newcommand*\flip[1]{ \overset{#1}{\leadsto}}

\newtheorem{theorem}{Theorem}
\newtheorem{lemma}{Lemma}
\newtheorem{corollary}{Corollary}
\newtheorem{claim}{Claim}

\newcommand{\qedclaim}{\hfill $\diamond$ \medskip}
\newenvironment{proofclaim}{\noindent{\em Proof of the claim.}}{\qedclaim}

\title{Exact number of flips required to sort a burnt stack of pancakes}

\author {
    Gerold Jäger\equalcontrib,
    Nacim Oijid\equalcontrib,
}
\affiliations {
Department of Mathematics and Mathematical Statistics, Umeå University, Umeå, Sweden \\
   gerold.jager@umu.se, nacim.oijid@umu.se
}

\usepackage{bibentry}

\begin{document}

\maketitle

\begin{abstract}
In this work, we consider the burnt pancake problem,
which is a well-studied  problem going back to 
a work of Gates and Papadimitriou from 1979. 
The problem is to sort a stack
of~$n$ one-sided burnt pancakes of different sizes,
by a sequence of flips
of the top pancakes,
such that at the end of the
flipping sequence the pancakes have increasing
size and the burnt sides
of all pancakes are face-down. The pancakes are
denoted  by $ 1,2,\dots,n$,
and a number is multiplied by $-1$, if the corresponding pancake has burnt side face-up.

Let $T(n)$ be the minimum number of flips to 
sort a special stack of $n$ 
pancakes, namely 
$\overline{I_n} := [-1,-2,...,-n]$. 
The instance $\overline{I_n} $
has strong relevance because of its easy structure and as it has been shown 
to be a worst-case instance for several small $n$.
Heydari and Sudborough gave in 1997 the currently best 
upper bound of $T(n)$, namely 
$ \lfloor (3n+3)/2 \rfloor $ for 
$ n \equiv 3 \bmod 4$, which later has been 
shown to be exact by a work of Cibulka from 2011. Except these two
works, no progress regarding lower and upper bounds has been made until now. 
In our work, we present that 
$ \lfloor (3n+3)/2 \rfloor $ 
is also an upper bound of $T(n)$ 
for $n \equiv 1 \bmod 4 $, which 
again matches the lower bound of Cibulka in 2011 
and thus is exact.
Furthermore, we show that our 
construction approach for 
$n \equiv 1 \bmod 4 $ and the one of Heydari
and Sudborough for $n \equiv 3 \bmod 4 $ 
cannot be applied for even $n$.
However, as there might be different 
construction approaches,
the case of even $n$ remains an open problem, 
where two possible values
for $T(n)$ are possible, namely 
$ (3/2)n  + 1 $ or 
$  (3/2)n  + 2 $. 
Finally, we found two values,
namely $n=24$, $n=26$,
where the lower bound is attained. 
\end{abstract}

\begin{links}
    \link{Code}{https://github.com/noijid/burnt-pancake}
\end{links}

\section{Introduction}

The problem of this work
can be motivated as follows.

For the buffet, the waiter of a restaurant 
gets a large stack of pancakes from the overworked cook. 
As usual, one side is burnt, and as the level of batter decreases, 
the pancakes became smaller and smaller.
Hence, the waiter ends up with a stack of one-sided burnt 
pancakes sorted by size, with the largest at the bottom and 
burnt side to be face-up. However, the waiter cannot 
serve them this way. He needs to turn all the burnt sides down, 
without changing the order. Having only a spatula, he can only 
perform flips to the top of the stack. How can he perform 
this transformation in a minimum number of flips?

The \emph{pancake sorting problem} was first mentioned by Dweighter (a pseudonym for J.E. Goodman) 
in~\cite{Dwe75} and its study was formally introduced by~\cite{GP79} under two variants: unburnt and burnt.
For a stack of $n$ pancakes of different sizes, the unburnt variant 
of this problem is to find the minimum number 
of flips required
to sort it, where each flip reverses the order of
the top $k$ pancakes for some $ 1 \le k \le n$. The burnt variant considers the extra condition
that each pancake has a burnt side and the stack needs both to be sorted and all 
the pancakes have their burnt side to be face-down. 
The (un)burnt $n$-pancake graph is defined as the graph 
containing one vertex per possible stack of $n$ pancakes, 
and two vertices are connected if it is possible to transform 
one stack to the other one in a single flip. 

Despite the simplicity of the problem, several applications have been found,
in binary tree simulation~\cite{AK86},
genome arrangements~\cite{FLRTV09} and parallel 
processing~\cite{HS97}. From a theoretical point of view, the study of this problem
is interesting both regarding structure and algorithmics. 
From a structural point of view, the pancake stack and its flipping sequences is a simple way to
represent large objects. 
The pancake graph
is of theoretical interest as it is a Cayley graph, 
and its structural properties have been an important area of research 
in structural graph theory~\cite{KANEVSKY1995, Qiu, Wang2024}. 
On the algorithmic side, the unburnt pancake sorting problem has been shown to be 
$\NP$-hard~\cite{BFR15}, whereas the complexity 
state of the burnt pancake sorting problem is open.
Furthermore, efficient algorithms have been developed to handle some restrictions 
or approximations of the problem~\cite{Seo2015}.

The investigations of the two variants of the problem has led to the 
study of two functions, $f(n)$ and $g(n)$ 
that count the minimum number of flips to sort any pancake stack 
in the unburnt and burnt variant, respectively.

When they introduced the problem formally, \cite{GP79}~proved that $\frac{17}{16} n \le f(n) \le \frac{5n+5}{3}$.
Whereas  Chitturi {\em et al.} improved the upper bound to  
$\frac{18}{11}n$~\cite{CFMMSSV09},  \cite{HS97} improved the lower bound to 
$\frac{15}{14} n $ and~\cite{peczarski2025} 
further by the additive constant $\frac{9}{14}$.
\cite{HS97}~also studied the worst-case instance regarding the lower bound
and showed an upper bound of $ \frac{9}{8}n+2$ for it. Recently, Amano~\cite{Ama24} improved 
this upper bound to $ \frac{15}{14}n+2$.

Considering the burnt pancake problem, \cite{GP79} 
showed that $\frac{3}{2} n-1 \le g(n) \le 2n+3$.
Let $T(n)$ be the minimum number of flips to sort 
the stack of pancakes $\overline{I_n}$,
which is the pancake stack in which all the pancakes are in increasing 
order with regard to their size, but with their burnt side face-up instead of
down. 
Naturally, it implies that $T(n) \le g(n)$. 
\cite{CB95} 
proved that, 
for $n \ge 10$, it holds that $\frac{3}{2} n \le T(n) \le  g(n) \le 2n-2$. 
They also conjectured that $\overline{I_n}$ is always a worst-case burnt pancake stack, 
i.e., one requiring the largest number of flip to be sorted.
They proved that this conjecture holds for $n \le 10$. 
The upper bound on $T(n)$ was improved 
by~\cite{HS97} who proved that 
$T(n) \le \frac{3n+3}{2}$ if $n \equiv 3 \bmod 4$, being very close to the lower bound. 
This upper bound appeared to be optimal when~\cite{Cib11} 
proved that it is also the general lower bound.
Hence, only the values for $n \equiv 0,1,2 \bmod 4$ remained to be determined.
\cite{Bou16} listed all values
$T(n)$ until $ n \le 27$ and 
upper bounds for $ n=28,29,30$, where
as we mention later, we have shown
that the values for 
$ n=24$ or $n=26$ are wrong, but all others are correct.

Finally, computing the values of $T(n)$ is challenging on its own, and we
refer the reader to the sequence~A078942 of the On-Line Encyclopedia of Integer Sequences~\cite{Oei25} to have a more general overview of the known values.

In the same article as mentioned above, \cite{Cib11}~disproved the conjecture that $\overline{I_n}$ 
requires the largest number of flips with a counterexample for $n = 15$. 
However, until now, $n = 15$ is the only known value for which $\overline{I_n}$ is 
not a stack requiring the largest number of flips, and all the worst-case
instances until $n =22$ have
been tested by~\cite{Bou16}. 

In this work, we show the same 
upper bound for $T(n)$ as in~\cite{HS97},
namely $ \frac{3n+3}{2} $, but for the case 
$ n \equiv 1 \bmod 4$. Together
with the case
$ n \equiv 3 \bmod 4$ from~\cite{HS97},
the whole case of odd $ n $ is solved.
Using results of~\cite{Cib11}, this also 
restricts the case of even $n$ to two
possible values,
namely $ \frac{3}{2}n+1$
or $ \frac{3}{2}n+2$.
For $n=24$ and $n=26$ we found
flipping sequences which attain
this lower bound, namely with $37$
and $40$ flips, respectively. 
Note that this disproves 
the values given for $n=24$ and $n=26$ in~\cite{Bou16}
and note that these two values are the only known
ones attaining this lower bound except the  value $n=2$.

Our main result is obtained
by an efficient algorithmic search, which
explores only few states of the pancake graph.
This kind of efficient search could find other applications
in problems whose encoding is much smaller than
the number of states to consider, as it is the case
in reconfiguration and enumeration problems.

This work is organized as follows.
First, we present
notations and basic results 
which are needed in our main
result and its proof. 
Then we present
our main result, a search algorithm
leading to this result and a sketch of its
proof. We also prove that the construction
of our main result does not work for
even $n$.
We close this work with some
conclusions and suggestions for future
work.

We provide an additional online
supplemental which includes the following:
\begin{itemize}
\item Three pdf files with full proofs of the
correctness of the flipping sequences, 
\item three pdf files
consisting of example instances for each of the investigated cases 
$ 1, 5, 9 \bmod 12$, namely $n=61,53,57$, 
in which each step of our constructed
flipping sequences is listed,
\item a
source code in Python which produces optimal flipping sequences for odd $n$
and optimal (or possibly optimal up to an additive factor $1$) for even $n$
explicitly,
\item
a source
code in Python which is able to check correctness of
the flipping sequences
for $ n \equiv 1 \bmod 4 $, i.e., of our main result 

(note that this correctness test
can show only correctness
of a given flipping sequence,
and note that this correctness
test works also for all other $n$, namely the case 
$ n \equiv 3 \bmod 4 $ which has been solved by~\cite{HS97} and the  case of even $n$).
\end{itemize}

    \section{Notations and basic results}
    \label{sec_not}
    
    \subsection{Notations}

    We will use the notations of~\cite{Cib11}, that we recall in the following.
    
    We denote by $I_n$ the sorted stack of pancakes, and by $\overline{I_n}$ the
    burnt pancake stack of size $n$ in which all pancakes are in increasing order with regard to their size, but
    burnt face-up. i.e., we have $\overline{I_n} = [-1, -2, \dots, -n]$.

    If $S = [x_1, x_2, \dots, x_s]$, where $ 1 \le s \le n$, is a stack of pancakes, we denote by $-S$
    its flip, and by $\overline{S}$ its opposite i.e., $-S = [- x_s, -x_{s-1},\dots, -x_1]$
    and $\overline{S} = [-x_1, -x_2,\dots, -x_s]$. 
    We denote by $T(n)$ the minimum number
    of flips required to sort $\overline{I_n}$.
    
In some parts, 
it will be helpful to add the (virtual)
pancake $\pm (n+1)$ which cannot be moved.

 An adjacency is a pair of consecutive pancakes that are already in
the right relative order: either both face-up and in increasing order $[s,
s+1]$, or both face-down and in decreasing order $[-(s+1), -s]$, for some $1
\leq s \leq n$. In other words, no flip is needed between them.
An {\em
anti-adjacency} is a pair of consecutive pancake in the same relative order with
respect to $\overline{I_n}$: either both face-up and in decreasing order $[s+1,
s]$, or both face-down and in increasing order $[-s, -(s+1)]$, for some $1 \le s
\le n $. In other words, no flip has been performed between them from the
original stack $\overline{I_n}$. In both adjacency and anti-adjacency, we also
include the virtual pancakes $\pm (n+1)$ to handle boundary cases.

    A set of consecutive pancakes is {\em a block} of size $s \ge 2$ if it is of the
    form $[x, x+1, x+2, \dots, x+s-1]$,
    where $x$ can be either positive or negative,
    and it is {\em a clan} of size $s \ge 2$ if
    it is of the form $[x, x-1, \dots, x - s +1]$. Note that clans and blocks are
    necessarily disjoint, and that $S$ is a clan if and only if $\overline{S}$ is a
    block. A pancake that is neither in a block nor in a clan is said to be {\em
    free}.
    Note that $\overline{I_n}$ is a clan of size $n$ and $I_n$ is a block of size
    $n$.

{\bf Example:} Consider the stack $[-5,-4,-1,-2,3,6,7,8]$. In this stack, $3$ is a free pancake, $-5,-4$ and $6,7,8,9$ are blocks (where the pancake $9$ is a virtual pancake added at the end of the stack to show that the pancake "$8$" is well placed) and $-1,-2$ is a clan.

    If $H_1$ and $H_2$ are two stacks, we denote by $H_1 \flip{k} H_2$ if $H_2$
    is obtained by one flip $(k)$ of $k$ pancakes in $H_1$. 

    A flip $H_1 \flip{k} H_2$ is said to be an {\em improve} if the number of 
		adjacencies in $H_2$ is strictly larger than it is in $H_1$. A flip that is not an improve is said to be a {\em waste}.

    \subsection{Basic results}
    
    The following results are used 
    to reach our main result. 
    
    We start with a result which corresponds to our main result,
    but holds for the case $ n \equiv 3 \bmod 4 $ (instead of 
    $ n \equiv 1 \bmod 4 $ in our
    result).
    
    \begin{theorem}{\cite[Observation 3]{HS97}} 
    \label{3mod4}
        Let $n \ge 23$ be an integer with $n \equiv 3 \bmod 4$. Then it holds that $T(n) = \frac{3}{2}n + \frac{3}{2}$.
    \end{theorem}

    It can easily be seen that
    each flipping sequence
    leading to $I_n$ needs 
    to have the flip of the entire
    stack
    at least twice.
    The following theorem 
    states that there is even a 
    flipping sequence which 
    begins with the flip of the entire stack.

    \begin{theorem}{\cite[Theorem~6]{CB95}}
    \label{start_with_full_flip}
        There exists a shortest sequence of flips for sorting $\overline{I_n}$,
        that begins with the flip of the entire stack.   
    \end{theorem}
    
    The following theorem states
    that the function $T(n)$
    increases not more than by $2$
    if $n$ increases by $1$. Having
    results for $T(n)$ 
    for the case of odd $n$,
    we will use this theorem to obtain also an upper bound for the case of even $n$. 
    
    \begin{theorem}{\cite[Theorem~7]{CB95}} \label{Tn+2}
        Let $ n\ge 1$ be an integer. Then it holds that $T(n+1) \le T(n) +2$. 
    \end{theorem}
    
    Note that the proof of 
    Theorem~\ref{Tn+2} is constructive, i.e., 
    given a flipping sequence for $\overline{I_n} $
    with $T(n)$ flips, where $ n\ge 1 $ is an integer, you can construct a
    flipping sequence for $\overline{I_{n+1}} $ with $ T(n)+2 $ flips.

    \section{Main result and method}
    \label{sec_main}
    
    Our main result is as follows.
    
    \begin{theorem}\label{theorem1mod4}
        Let $n \ge 29$ be an integer with $n \equiv 1 \bmod 4$. Then it holds that $T(n) = \frac{3}{2}n + \frac{3}{2}$.
    \end{theorem}

    We prove this result
    in the additional online supplemental.
    However, we have independently tested its
    correctness by a computer
    program, also provided in the additional online supplemental, until $n=1501$.

    In this main result,  
    exact values for the specific instance
    $\overline{I_n}$ of the burnt pancake sorting problem
    are computed, namely for $ n \equiv 1 \bmod 4 $. On the one hand, the proof of
    this theorem 
    is based on a precise description of the flipping sequences together with the states of
    the stack in each step. However, these flipping sequences required an advanced algorithm to be
    obtained. We will present this algorithm in the next section.

We conducted our computational experiments on a Linux-based high-performance computing cluster. Each compute node 
was equipped either with two Intel® Xeon® Gold 6132 processors running at 2.6 GHz or two AMD EPYC 9754 (Zen 4) processors running at 2.25 GHz. For each run, we used 16 physical CPU cores and 64 GB of RAM. Our experiments consisted of 20 independent runs executed in parallel.
    
    \subsection{Algorithmic search of the sequences}
    
    This search can be divided into two steps. The first one is to reach optimal flipping sequences for burnt stack of pancakes of sizes
    $29, 33$ and $37$, and the second one to generalize these sequences to sizes $ 5 \bmod 12, 9 \bmod 12$ and $1 
    \bmod 12 $.

 The main idea behind the use of 
$"\!\!\bmod 12"$ is that the flipping
sequence for $n = 3$ is $[3,2,3,2,3,2]$ and the one for $n = 15$ is 
$[15,10,4,6,14,6,4,10,15,10,4, 6,$ $14,6,4,10,15,10,4,6,14,6,4,10]$ (see~\cite{HS97}), which is
obtained from the previous one by either adding 12 (the number of new pancakes)
to the value or by inserting $[10,4,6]$ or $[6,4,10]$. This is the first value for
which we obtained this simplicity in the computation from a value to the next
one (modulo 12), which made it easy to implement and testing if this phenomenon 
can appear in other flipping sequences. It turns out that this is the case.
    
While
analyzing the flipping sequences of~\cite{HS97}, for $n = 3 \bmod
4$, we realized that the main arguments making this flipping
sequence sort $\overline{I_n}$ are based on the assumption that the number of
pancakes is odd. Thus, it was natural to try to extend their result from $3 \bmod
4$ to $1 \bmod 4$. However, the known values on $T(n)$ showed that, if we want to
{prove that $T(n) = \frac{3}{2}n +\frac{3}{2}$ for $n \equiv 1 \bmod 4$, it would be necessary
to start at least at $n = 21$, since $T(17) = 28 = \frac{3}{2}n + \frac{3}{2} + 1$ (see~\cite{CB95,Bou16}).

\subsubsection*{Base case}

The first step was to compute the base cases that can be generalized to larger
values. In order to have a generalizable flipping sequence,
following a similar idea to the flipping sequence of~\cite{HS97}, we aimed
to find a flipping sequence that first starts with all the wastes, and after that only
performs improves. 

We first restricted ourselves to the flipping sequences that start with a flip of
the entire stack. Using  
Theorem~\ref{start_with_full_flip},
we know that such a flipping sequence exists. 
The next wastes cannot be performed in any order, but we have the following
lemma that implies some structure on any flipping sequence that starts with all
the wastes.
\begin{lemma}\label{clan of size 3}
    Let $P$ be a stack of pancakes containing a clan $C$ as a substack of size $|C| \ge 3$. Then $P$ cannot be sorted without waste.
\end{lemma}

\begin{proof}
    Let $x$ be an integer and 
    \\
    $ C = [x, x-1, x-2]$ \big/
    $ [-(x-2), -(x-1), -x] $ 
    
    be a clan of size $3$ in $P$. If
    the first flip modifying $C$ is a waste, then there is nothing to do.

    Otherwise, the resulting stack is either 
    \\
    $\big[-(x-1), -x, \dots, x-3,
    x-2 \big]$ or  $\big[x-1, x -2, \dots, -(x+1), -x \big]$. 
    In both cases,
    the next flip cannot be an improve since to be an improve the pancake that follows $-(x-1)$ ($(x-1)$ resp.) after the next flip must be $x$ ($-(x-2)$ resp.), which currently has
    the wrong side up in the stack. 
\end{proof}

In total, we aim at having
$\frac{3}{2}n  + \frac{3}{2}- n = \frac{n}{2} +
\frac{3}{2} $ wastes, 
and because of Lemma~\ref{clan of size 3}, 
we have to perform them  
before reaching the state where all clans have size at most~$2$. 
Since one waste was
already used to flip the entire stack, and since we need at least
$\frac{n-1}{2}$ wastes to split the stack into clans of size~$2$, this forces
all the wastes but one to be used to split clans. Following the idea 
of~\cite{HS97}, we again restricted ourselves to the flipping sequences that
have a flip of the entire stack as last waste. If such a flipping sequence
exists, this last waste ensures a very precise structure of the rest of the
flipping sequence as it will be described later. 

After all these hypotheses, we know that we can generate candidate flipping
sequences as follows, for a permutation $\sigma$ of size $ \frac{n}{2} - \frac{3}{2}$.
\begin{itemize}
    \item We flip the entire stack.
    \item For $1 \le i \le \frac{n}{2} - \frac{3}{2}$, we flip between the clan $[2\sigma(i)+1, 2\sigma(i) ]$ and the clan $[2\sigma(i)+3, 2\sigma(i) +2]$
    \item We flip the entire stack.
\end{itemize}
Note that this algorithm can easily be parallelized since each permutation has to be handled separately.

We implemented this algorithm in Python
and run the corresponding
computer program for $n = 21$ and $n = 25$, but we did not obtain the desired flipping sequences. Hence, our results 
started at $n = 29$.
For $n = 29$, this leads us to $13! \simeq 6.2 \cdot 10^{9}$ 
flipping sequences
to try. For each of them, since we know that no more wastes are allowed, we can
easily check in quadratic time if it is possible to sort the resulting 
stack of pancakes using only improves:
each time a flip is an improve, if we denote by $i$ the first pancake of
the stack, we know that the flip has to occur before pancake $(- i )+ 1$. 
This computer program provided us $20$ flipping sequences of size $45$, which were obtained
in $31$ hours using $16$ cores. Among them, we chose one arbitrarily and we
managed to generalize it to any values $n \equiv 5 \bmod 12$ as we will describe in the next
subsection.

For $n = 33$ and $n = 37$, this would lead us to $15! \simeq 1.3 \cdot 10^{12}$
and $17! \simeq 3.6 \cdot 10^{14}$ permutations to test. Compared to the values
for $29$ the computation time would not end in reasonable time. To solve this
problem, we randomized the generation of the permutations, hoping that the
number of correct flipping sequences will also increase sufficiently so that we can
find one with this random search. It happens that we reached both 
flipping sequences for $33$ and $37$ in less than two days of computation. 

\subsubsection*{Generalized sequences}

Generalizing the flipping sequences from $n = 29$, $33$ and $37$ to any $n \equiv 1
\bmod 4$ was
mostly done inductively
with the help of another
computer program in Python. For sake of simplicity, we explain here the process
that generalizes the flipping sequence from $n = 29$ to all  larger values $n \equiv5 \bmod 12$. 
We have applied the same process from $n = 33$ to all larger values  $ n \equiv 9 \bmod 12$ and
from $ n = 37$ to all larger values $n \equiv 1 \bmod 12$. 

The main idea is to consider the wastes of the
flipping sequences for $n = 15$ which are $[15, 10, 4, 6, 14, 6, 4, 10,
15]$ and to observe that a good solution to go from $29$ to $41$ is to
introduce 
the subsequence of three wastes $[10, 4, 6]$ 
somewhere at the beginning, and the subsequence of three wastes $[6, 4, 10]$
somewhere at the end of the waste sequence for $n=29$. 
The intuition is that the flipping sequence for $n=15$ is obtained by adding
these wastes to the sorting sequence for $n = 3$, hence we can naturally think
that these wastes placed at the correct place in the flipping sequence manage to
handle the $12$ new pancakes using exactly~$6$ wastes and then we have $12$ 
additional improves.

Denote by $W$ the set of wastes of a flipping
sequence. 
The idea is that 
inserting these two subsequences into $W$ will split it into
three (possibly empty) sets $W_1$, $W_4$ and $W_6$. Then, $W_2 = [10, 4, 6]$ is inserted
between $W_1$ and $W_4$, and $W_5 = [6, 4, 10]$ is inserted between $W_4$ and $W_6$
($W_3$ will be defined later as we will need one more observation to introduce
it, see the end of this paragraph). 
The computation has shown that to apply this method a second time, i.e., from $n=41$ to $n=53$, works naturally for $W_5$, which
is always inserted in the same place going from any integer $k$ to the integer
$k+12$. I.e., $W_5$ for $n = 53$ is obtained by inserting a
new $[6,4,10]$ subsequence before $W_6$ and adding $12$ to the flip $(10)$ of the
previous $[6,4,10]$ subsequence in the flipping sequence of $41$, i.e., adding $[6,4,22,6,4,10]$ before $W_6$ in the flipping
sequence of $n=29$.
However, the next insertions of $W_2$ will happen in the previous $[10,4,6]$
subsequence by adding also $12$ to all their values. I.e., the part
$[10,4,6]$ of the flipping sequence for $n=41$ becomes $[10+12, 4 +12, 10, 4,
                                                         6, 6+12] = [22, 16, 10, 4, 6, 18]$ in the flipping sequence for $53$. Hence, to
simplify the notation, we divide $W_2$ into two parts, and we now set $W_2 =
[10,4]$ and $W_3 = [6]$, such that the flipping sequence is easier to describe.

However, even by knowing which subsequence of flips to insert, finding their positions
in the flipping sequence is not so easy. Hence, we used the computer program to do
so.
This aims to find exactly how many flips should be in $W_1$ before inserting $W_2$ and $W_3$, and how many flips should be in $W_4$ before inserting $W_5$.
All other values of the
sequence should not be modified, but, depending on whether we count from the
beginning or the end of the stack, some flips may require a ``$+12$'' in their value.
Hence, we tested all the possibilities to obtain the general form of $W_1$ and $W_4$.
We generated all the possible flipping sequences that satisfy all these
properties, and tested whether it was possible to complete them with only improves.
Note that this algorithm is much faster than the  program to find a flipping sequence for $n=29$ 
mentioned in the last section, as it only
requires to test $(\frac{n-1}{2})^2$ positions to insert the subsequences
$[10,4,6] $
and $[6,4,10]$
(which leads to the positions of $W_2, W_3$ and $W_5$ in the flipping sequence). Moreover, while going from $n = 29$ to $ 41$, each value of $W_1, W_4$ and $W_6$ is preserved or increased by $12$. Thus, only $2^{\frac{n-1}{2} -1}$ flips have to be tested for these values. 
Once one flipping sequence was found
this way, we tried by hand to generalize it by repeating the same process to
reach the general formula. All this process is possible, since testing whether a
flipping sequence is correct can be done in quadratic time.

\subsection{Flipping sequences} \label{flipping sequence}

In the following we provide optimal flipping sequences to sort the stack $\overline{I_n}$
for any $n \ge 29$ with $n \equiv 1 \bmod 4$. Similarly to the
flipping sequences of~\cite{HS97}, our flipping sequences depend on the
value of $n \bmod 12$. Due to the length of the proof and as they are mostly
based on a technical analysis of the stack structure after each flip, we only
provide the flipping sequences and an overview of how they are built. We refer the reader to the additional online supplemental to have the full versions of the
proofs, with examples of optimal flipping sequences applied to stacks of
$61$, $53$ and $57$ pancakes, and
with source codes providing the flipping sequences and checking
their correctness.

\subsubsection{The case $n \equiv 1 \bmod 12 $} \hspace{.1cm}

Let $n$ be an integer such that $n \equiv 1 \bmod 12$ and $n \ge 37$ and let $s =
\frac{n-37}{12}$. We define the flipping sequence $S_1 = W^1 + A^1 + B^1$ as
follows:
\begin{itemize}
    \item $W^1 = W^1_1 + W^1_2 +  W^1_3 + W^1_4 + W^1_5 + W^1_6$:
\begin{itemize}
    \item $W^1_1 = \{n, n-9, n-21, n-13, n-15\}$,
    \item $W^1_2 = \{n-39-6i\}_{0 \le i \le 2s-1}$,
    \item $W^1_3 = \{6+12i\}_{0 \le i \le s-1}$,
    \item $W^1_4 = \{n-19, 10, n-1, n-13, n-15, 14, n-9, n-19, n-15, 12, n-3, 2, 14, n-17\}$,
    \item $W^1_5 = \{6,4, n-39-12i\}_{0 \le i \le s-1}$,
    \item $W^1_6 = \{n\}$.
\end{itemize}

    \item $A^1 = A^1_1 + A^1_2 + A^1_3 + A^1_4 + A^1_5$:
\begin{itemize}
    \item $A^1_1 = \{6, 32, 14, 20, 10, 2, 28, 4, 30, 8\}$,
    \item $A^1_2 = \{44 + 12i, 4, 6\}_{0 \le i \le s-1}$,
    \item $A^1_3 = \{n-1\}$,
    \item $A^1_4 = \{n-9-6i, n-11-12i, n-5-6i\}_{0 \le i \le s-1}$,
    \item $A^1_5 = \{\frac{n-1}{2}, \frac{n-1}{2}+4, \frac{n-1}{2}-2, 10, \frac{n-1}{2}-4, \frac{n-1}{2}+2, \frac{n-1}{2}+8, n \}$.
\end{itemize}

    \item $B^1 = B^1_1 + B^1_2 + B^1_3$:
\begin{itemize}
    \item $B^1_1 = \{\frac{2n+1}{3}+3\}$,
    \item $B^1_2 = \{\frac{2n+1}{3} -3 -8i, \frac{2n+1}{3} -1 -8i, \frac{2n+1}{3} + 7 +4i, 6+12i, 4+12i, \frac{2n+1}{3}+3 +4i\}_{0 \le i \le s-1}$,
    \item $B^1_3 = \{10, 24, 2, n-3, n-15, n-19, n-17, n-11, n-27, n-9, n-1, n-21,$ $ n-9, 4, 6, 18, 14\}$.
\end{itemize}

\end{itemize}

\begin{lemma}\label{lemma 1 mod 12}
    Let $n \ge 37$ be an integer with $n \equiv 1 \bmod 12$. The sequence $S_1$ sorts $\overline{I_n}$ and has size $\frac{3}{2}n+
    \frac{3}{2}$. 
    \end{lemma}

\subsubsection{The case $ n \equiv 5 \bmod 12 $} \hspace{.1cm}

Let $n$ be an integer such that $n \equiv 5 \bmod 12$ and $n \ge 29$. Let $s = \frac{n-29}{12}$. We define the flipping Sequence $S_5 = W^5 + A^5 + B^5$ as follows:

\begin{itemize}
    \item $W^5 =  W^5_1 + W^5_2 + W^5_3 + W^5_4 + W^5_5 + W^5_6$:
\begin{itemize}
    \item $W^5_1 = \{n, n-9, 4, 10, n-3, n-15, n-25\}$,
    \item $W^5_2 = \{n-31-6i\}_{0 \le i \le 2s-1}$,
    \item $W^5_3 = \{6 + 12i\}_{0 \le i \le s-1}$,
    \item $W^5_4 = \{n-7, n-1, n-11, 6, 4, n-5, 10, n-15\}$,
    \item $W^5_5 = \{ 6, 4, n - 31 - 12i\}_{0 \le i \le s-1}$,
    \item $W^5_6 = \{n\}$.
\end{itemize}
    
    \item $A^5 = A^5_1 + A^5_2 + A^5_3 + A^5_4 + A^5_5$:
 \begin{itemize}
    \item $A^5_1 = \{18, 20, 12, 4, 6, 24, 16, 26, 16, 4, 18\}$,
    \item $A^5_2 = \{36 +12i, 4, 6\}_{0\le i \le s-1} $,
    \item $A^5_3 = \{n-1\}$,
    \item $A^5_4 = \{n-9-6i, n-11-12i, n-5-6i\}_{0\le i \le s-1}$,
    \item $A^5_5 = \{\frac{n-1}{2} -6, \frac{n-1}{2} + 6, n\}$.
\end{itemize}
    \item $B^5 = B^5_1 + B^5_2 + B^5_3$:
 \begin{itemize}
    \item $B^5_1 = \{\frac{2n-1}{3} -3\}$,
    \item $B^5_2 = \{\frac{2n-1}{3} -9 - 8i, \frac{2n-1}{3} - 7 -8i,\frac{2n-1}{3}+1 +4i, 6 +12i, 4+12i, \frac{2n-1}{3} -3 + 4i\}_{0 \le i \le s-1}$,
    \item $B^5_3 = \{12, n - 5, 4, n - 13, n-15, n-21, n - 25, n-11, n-9, n-1, 6, 4, 14\}$.
\end{itemize}
\end{itemize}

\begin{lemma}\label{lemma 5 mod 12}
    Let $n \ge 29$ be an integer with $n \equiv 5 \bmod 12$. The sequence $S_5$ sorts $\overline{I_n}$ and has size $\frac{3}{2}n + \frac{3}{2}$.
\end{lemma}

\subsubsection{The case $n \equiv 9 \bmod 12 $} \hspace{.1cm}

Let $n$ be an integer such that $n \equiv 9 \bmod 12$ and $n \ge 33$. Let $s = \frac{n-33}{12}$. We define the flipping Sequence $S_9 = W^9 + A^9 + B^9$ as follows:
\begin{itemize}
 \item $W^9 = W^9_1 + W^9_2 + W^9_3 + W^9_4 + W^9_5 + W^9_6$:
\begin{itemize}
    \item $W^9_1 = \{n, 14, 4, 10, n-7, n-29\}$,
    \item $W^9_2 = \{n - 35 -6i\}_{0 \le i \le 2s-1}$,
    \item $W^9_3 = \{ 6 +12i\}_{0 \le i \le s-1}$,
    \item $W^9_4 = \{n-11, n-1, n-9, 4, n-13, n-11, 8, 10, n-5, 10, n-19\}$,
    \item $W^9_5 = \{6, 4, n-35-12i\}_{0 \le i \le s-1}$,
    \item $W^9_6 = \{n\}$.
\end{itemize}

    \item $A^9 = A^9_1 + A^9_2 + A^9_3 + A^9_4 + A^9_5$:

\begin{itemize}
    \item $A^9_1 = \{24, 16, 14, 4\}$,
    \item $A^9_2 = \{40 +12i, 4, 6\}_{0 \le i \le s-1}$,
    \item $A^9_3 = \{n-1\}$,
    \item $A^9_4 = \{n - 9 -6i, n - 11 - 12i, n-5 - 6i\}_{0 \le i \le s-1}$,
    \item $A^9_5 = \{\frac{n-1}{2} +10, 12, 10, 4, 14, 22, 6, 24, \frac{n-1}{2}+14, \frac{n-1}{2}-6, \frac{n-1}{2}, n\}$.
\end{itemize}

    \item $B^9 = B^9_1 + B^9_2 + B^9_3$:

\begin{itemize}
    \item $B^9_1 = \{\frac{2n}{3}\}$,
    \item $B^9_2 = \{\frac{2n}{3} -6 - 8i, \frac{2n}{3} -4 -8i, \frac{2n}{3} +4 +4i, 6 +12i, 4 +12i, \frac{2n}{3} +4i\}_{0 \le i \le s-1}$,
    \item $B^9_3 =  \{18, n-5, n -19, n-13, 8, n-19, n-29, n-9, 14, n-1, n-17, n-13, 10, 4, n-9\}$.
\end{itemize}

\end{itemize}

\begin{lemma}\label{lemma 9 mod 12}
    Let $n \ge 33$ be an integer with $n \equiv 9 \bmod 12$. The sequence $S_9$ sorts $\overline{I_n}$ and has size $\frac{3}{2}n+\frac{3}{2}$.
\end{lemma}

\subsection{Description of the flipping sequences} \label{subsec: overview}

Even though the flipping sequences are different depending on the value of $n \bmod 12$,
they are divided in three parts, $W, A$ and $B$, where each of them is divided in
several steps. 

Part $W$ contains $ \frac{n+3}{2}$ flips.  It aims first to split the stack into $\frac{n-1}{2}$ clans of size $2$,
and the free pancake $1$. Then it performs a full flip of the stack to put the
free pancake $-1$ to the top of the stack. To achieve this result, Step $W_1$
and Step $W_4$ are the ones that create the clans for a {\em small} stack, i.e.,
of size $29, 33$ or $37$, respectively. They are generalized to handle the
first and last pancakes of the stack, while for 
large~$n$, most of the pancakes  are handled by the other steps. When the size of the stack increases,
Step $W_2$ creates new clans of size~$6$, 
one set of clans starting from the larger pancakes (of value almost $n$) and one set of clans from the smaller ones (of constant value), until reaching pancakes of value $\frac{n}{2} +c$ for some constant $c$.
Step $W_3$ and Step $W_5$ aim to separate these clans of size $6$ into clans of size $2$. First,
Step $W_3$ handles the clans of one set to transform them into
a clan of size $2$ and one of size $4$. Then, Step $W_5$ handles the remaining clans of size $4$ and $6$, and separates them into regularly separated clans of size~$2$. Intuitively, the clans are split like they are in the flipping sequence for $n=15$.
Finally, Step $W_6$ consists of fully
flipping the stack to put the free pancake $-1$ to the top. This makes it
possible to only have improves after this flip.

 Then, part $A$ consists of $\frac{n+1}{2}$ improves to transform all the clans
 into blocks of size $2$. 
In particular, we start with the free pancake $-1$ to the top of the
 stack. Then, 
 for $ 1 \le i \le \frac{n+1}{2} $, the $i$-th flip of part $A$ starts with the pancake $-(2i-1)$ on
 the top of the stack and is flipped between the clan $(2i+1, 2i)$, transforming
 it into the block $(2i-1, 2i)$ and putting the free pancake $-(2i+1)$ to the top.
 The division of part $A$ in several steps is due to the structure of the stack
 achieved during Part $W$.

 Finally, part $B$ performs $\frac{n-1}{2}$ improves to merge these blocks
 together to obtain $I_n$. The flip of Step $B_1$ puts back this pancake of value $\frac{n}{2} +c$ 
 to the top of
 the stack. Then, Step $B_2$ merges all the blocks
 except the beginning and the end of the stack (which both correspond to a small
 stack of pancakes), ending up with a large block containing all the pancakes, except
 some of the smaller ones and some of the larger ones. 
 Step $B_3$ finally
 merges the remaining blocks that correspond to a small base stack, i.e., corresponding to a stack of size $29, 33$ or $37$,
 with the big block
 created during Step $B_2$.

\subsection{Consequences} \label{subsec: consequences}

\begin{corollary}\label{theorem1mod2}
    Let $n \ge 19$ be an odd integer. Then it holds that $T(n) = \frac{3}{2}n + \frac{3}{2}$.
\end{corollary}

\begin{proof}
    Let $n \ge 19$ be an odd integer. The lower bound is already known from~\cite{Cib11}, i.e., $T(n) \ge  \frac{3}{2}n + \frac{3}{2}$.

    For the upper bound, if $n=19,21,23,25,27$, the result is provided  by~\cite{Bou16},
    and we obtained it also independently by a Python program
    available in our online
    supplemental.
    If $n \ge 29$, depending on the value of $n \bmod 12$, the result is provided
    by Theorem~\ref{theorem1mod4}
    or by Theorem~\ref{3mod4}
    from~\cite{HS97}.
\end{proof}

\begin{corollary}\label{0 mod 2}
    Let $n \ge 14$ be an even integer. 
    We have: $$T(n) \in \left\{\frac{3}{2}n +1, \frac{3}{2}n +2\right\}.$$
\end{corollary}

\begin{proof}
Let $n\ge 14$ be an even integer.
    For $n=14,16,18$, the result is provided by~\cite{Bou16}, 
    and we obtained it also independently by a Python program
    available in our online
    supplemental.
    If $n \ge 20$, 
    by Theorem~\ref{Tn+2}, we know that $T(n+1) -2 \le T(n) \le T(n-1) +2$.
    Since $n$ is even, both $n-1$ and $n+1$ are odd. Thus, by
    Corollary~\ref{theorem1mod2}, we have $T(n+1) = \frac{3 }{2}(n+1) + \frac{3}{2} =
    \frac{3}{2}n + 3$ and $T(n-1) = \frac{3}{2}(n-1) + \frac{3}{2} = \frac{3}{2}n$, which
    leads us to:

$$ \frac{3}{2}n + 1 \le T(n) \le \frac{3}{2}n +2. \qedhere $$
\end{proof}

\section{Limit of the method} \label{sec:limit}

In this section, we prove that the method which we and~\cite{HS97} used to obtain flipping sequences 
for odd $n$ cannot be used to obtain a flipping sequence for even $n$ 
and to reach the best known lower bound, i.e., $\frac{3}{2}n +1$ 
flips.  However, even though these sequences do not exist for small values, 
by a computer program in C++,
based on BnB,
we found the following flipping sequences that sort $\overline{I_{24}}$ and
$\overline{I_{26}}$ and attain the lower bound, i.e., they have 37 and 40 flips,
respectively. Note that this disproves two values of Table 1
of~\cite{Bou16}.

\smallskip

\textbf{The flipping sequence for $\overline{I_{24}}$:} 24, 10, 12, 20, 5, 15, 5, 8, 19, 17, 5, 10, 14, 7, 22, 17, 6, 13, 10, 21, 23, 18, 11, 20, 6, 18, 24, 9, 7, 19, 5, 16, 23, 9, 15, 11, 18.

\textbf{The flipping sequence for $\overline{I_{26}}$:} 26, 12, 8, 25, 10, 14, 26, 21, 18, 13, 8, 24, 9, 12, 8, 22, 4, 7, 23, 19, 25, 11, 23, 16, 6, 11, 13, 9, 6, 22, 10, 14, 12, 18, 15, 24, 7, 14, 18, 21. 

\smallskip

For proving the mentioned
limitation,
we first need the following lemma that is not explicitly stated in this form
in~\cite{Cib11} but which is shown in the proof of his Lemma~4 as follows.

\begin{lemma} \label{cibulka3n2}
    Let $n$ be an even positive integer. If a flipping sequence of size $\frac{3}{2}n +1$ sorts $\overline{I_n}$, all its flips either create an adjacency or break an anti-adjacency.
\end{lemma}

To make the paper self-contained, we provide a short proof of this lemma
even though most of the arguments are from~\cite{Cib11}.

\begin{proof}
    Let $C$ be a stack of pancakes
    and let $C'$ be a stack obtained from a single flip from $C$.

    A block (clan resp.) is called a \emph{surface block (clan resp.)} if the topmost pancake of the stack is part of it, otherwise it is called \emph{deep}.\\
For a stack $C$ of pancakes, we define the value $v(C)$ as 

follows:~\footnote{Note that in comparison to~\cite{Cib11} for a stack of pancakes
$C$, we have interchanged the notations of $\bar{C}$ and $-C$.}

\begin{align*}
    v(C) \overset{\mathrm{def}}{=}
&a(C)-a^{-}(C) +l(C)-l^{-}(C)
-\frac13\bigl(b(C)-b^{-}(C)\bigr) \\
+ &\frac13\bigl(o(C)-o^{-}(C)\bigr)
+\frac13\bigl(ll(C)-ll^{-}(C)\bigr),
\end{align*}

where

\[
\begin{aligned}
a(C)
&\overset{\mathrm{def}}{=}
\text{number of adjacencies}
\\
l(C)
&\overset{\mathrm{def}}{=}
\begin{cases}
1 & \text{if the lowest pancake is $n$} \\
0 & \text{otherwise}
\end{cases}
\\
b(C)
&\overset{\mathrm{def}}{=}
\text{number of deep blocks}
\\
o(C)
&\overset{\mathrm{def}}{=}
\begin{cases}
1 & \text{if the pancake on top of the stack is the free $-1$} \\
  & \text{or if $1$ is in a block (necessarily with $2$)} \\
0 & \text{otherwise}
\end{cases}
\\
ll(C)
&\overset{\mathrm{def}}{=}
\begin{cases}
1 & \text{if the lowest pancake is $n$}\\
 & \text{ and the second lowest is $n-1$} \\
0 & \text{otherwise}
\end{cases}
\\
a^{-}(C)
&\overset{\mathrm{def}}{=}
a(\overline{C})
=
\text{number of anti-adjacencies in } C
\\
l^{-}(C)
&\overset{\mathrm{def}}{=}
l(\overline{C})
\\
b^{-}(C)
&\overset{\mathrm{def}}{=}
b(\overline{C})
=
\text{number of deep clans in } C
\\
o^{-}(C)
&\overset{\mathrm{def}}{=}
o(\overline{C})
\\
ll^{-}(C)
&\overset{\mathrm{def}}{=}
ll(\overline{C}).
\end{aligned}
\]

Cibulka~\cite{Cib11} proved that any flip changes the value of $v(C)$ by at most $4/3$, and that $v(I_n) = -v(\overline{I_n}) = n+\frac 23$. Hence, as $\frac 4 3(\frac  3 2 n +1) = 2n + \frac 4 3 = v(I_n) - v(\overline{I_n})$, a flipping sequence of size $\frac 3 2 n +1$ sorts $\overline{I_n}$ if and only if all the flips of the sequence increase the value of $v$ by $\frac 4 3$. We also define the following value to handle the different changes after a flip:

\[
\begin{aligned}
\Delta a
&\overset{\mathrm{def}}{=}
a(C')-a(C)
&
\Delta a^{-}
&\overset{\mathrm{def}}{=}
-\bigl(a^{-}(C')-a^{-}(C)\bigr)
\\
\Delta l
&\overset{\mathrm{def}}{=}
l(C')-l(C)
&
\Delta l^{-}
&\overset{\mathrm{def}}{=}
-\bigl(l^{-}(C')-l^{-}(C)\bigr)
\\
\Delta b
&\overset{\mathrm{def}}{=}
-\frac13\bigl(b(C')-b(C)\bigr)
&
\Delta b^{-}
&\overset{\mathrm{def}}{=}
\frac13\bigl(b^{-}(C')-b^{-}(C)\bigr)
\\
\Delta o
&\overset{\mathrm{def}}{=}
\frac13\bigl(o(C')-o(C)\bigr)
&
\Delta o^{-}
&\overset{\mathrm{def}}{=}
-\frac13\bigl(o^{-}(C')-o^{-}(C)\bigr)
\\
\Delta ll
&\overset{\mathrm{def}}{=}
\frac13\bigl(ll(C')-ll(C)\bigr)
&
\Delta ll^{-}
&\overset{\mathrm{def}}{=}
-\frac13\bigl(ll^{-}(C')-ll^{-}(C)\bigr).
\end{aligned}
\]

Cibulka proved that $\Delta a, \Delta a^-, \Delta l, \Delta l^- \in \{-1,0,1\}$ and $\Delta b, \Delta b^-, \Delta o, \Delta o^-, \Delta ll, \Delta ll^- \in \{-1/3, 0, 1/3\}$.

By definition of $a, a^-, l$ and $l^-$, if one of these four values increases, the flip either creates an adjacency or removes an anti-adjacency. Cibulka also proved the following implications:
\begin{align*}
    \Delta l = \Delta l^- =0 ~\&~  \Delta ll > 0 &\implies \text{it creates an adjacency} \\
    \Delta l = \Delta l^- =0 ~\&~ \Delta ll^- > 0 &\implies \text{it breaks an anti-adjacency}\\
    \Delta o >0 ~\&~ \Delta o^- >0 &\implies \Delta b = \Delta b^- = 0 
\end{align*}

Hence, if a flip does neither break an anti-adjacency nor create an adjacency, it satisfies that $\Delta a = \Delta a ^- = \Delta l = \Delta l^- = 0$ and at least three terms among $\Delta ll, \Delta ll^-, \Delta o, \Delta o^-, \Delta b, \Delta b^-$ are non-positive. Thus, this flips verifies that $v(C') - v(C) \le 1$, which cannot happen in a sorting sequence of size $\frac 3 2 n +1$.
\end{proof}

We now can prove that no flipping sequence with the structure we used can sort 
$-I_n$ in $\frac{3}{2}n+1$ flips when $n$ is even. 
This structure means starting with the flip $(n)$, 
then containing all the wastes, another flip $(n)$ 
and finally only improves. 
Note that in the odd case the second flip $(n)$ is not an improve, but instead the third one is.
(Note that in the even case, the third one
may not exist.)
Suppose by contradiction that 
such a flipping sequence exists.

\begin{claim}\label{claim:secondn-flip-improve}
    The second flip $(n)$ has to be an improve.
\end{claim}

\begin{proofclaim}
		 By Lemma~\ref{cibulka3n2}, if a flipping sequence of size $\frac{3}{2}n+1$ 
		 exists, all its flips have either to create an adjacency, or break an 
		 anti-adjacency. Hence, since the first flip $(n)$ 
		 already breaks the 
		 anti-adjacency between $-n$ and the end of the stack (a virtual 
         $-(n+1)$-pancake),  when the second flip $(n)$ is  performed, 
         the last pancake of the stack 
         cannot be the pancake $-n$. Thus, it cannot break an 
		 anti-adjacency, and therefore it has to create an adjacency
         with the virtual pancake $+(n+1)$, which is only 
		 possible if it puts the pancake $n$ to the end of the stack.
\end{proofclaim}

\begin{claim}\label{claim:breakn-n-1}
    The anti-adjacency $(n, n-1)$ has to be broken by a waste.
\end{claim}

\begin{proofclaim}
		Since the second flip $(n)$ puts the 
        pancake~$n$ to the end of the stack, the 
		flip before has to put the pancake $-n$ to its top. Suppose that 
		this flip is a flip $(k)$, for some integer $1 \le k \le n-1$. Before this 
		flip $(k)$, the pancake $n$ has to be positively in position $k$. Hence, 
		either the clan $(n, n-1)$ has already been broken, and there is nothing to 
		do, or the $(k+1)$-th pancake of the stack is the pancake $n-1$, and this 
		flip $(k)$ breaks this clan, which concludes the proof. In both cases, as this 
		happens before the second flip $(n)$, by hypothesis, the flip that breaks 
		this clan is a waste.
\end{proofclaim}

\begin{claim}\label{claim:break2-1}
    The anti-adjacency $(2,1)$ has to be broken by a waste.    
\end{claim}

\begin{proofclaim}
		Suppose by contradiction that it does not. As after the first flip $(n)$, the 
		two last pancakes of the stack are $[2,1]$ in that order, all the other 
		flips before the second flip $(n)$ have to be flips $(k)$ with $1 \le k \le n-2$. 
		Hence, the two pancakes $1$ and $2$ cannot be moved before the second flip  
		$(n)$, which will put $[-1,-2]$ on the top of the stack. However, by hypothesis, all 
		the flips after the second flip $(n)$ have to be improves, which is not 
		possible as the only pancake that can go after an improve when $-1$ is the 
		top pancake of the stack is the pancake $2$ which has the wrong side up in 
		the stack.
\end{proofclaim}

Since we supposed that a $(\frac{3}{2}n +1)$-flipping sequence exists, note that $n$ flips have to be improves, and we already determined three of the wastes: the first 
flip $(n)$, the ones that breaks the anti-adjacency $(n,n-1)$ and $(2,1)$. Hence $\frac{n}{2}-2$ wastes remain to be determined. 

\begin{claim}\label{claim:breakclan2k-2k-1}
    Let $1\le k \le \frac{n}{2}$ be an integer. There is one waste that breaks the clan $(2k, 2k-1)$.
\end{claim}

\begin{proofclaim}
    We already proved through Claim~\ref{claim:breakn-n-1} and Claim~\ref{claim:break2-1} that for $k=1$ and $k = \frac{n}{2}$, there must be a waste breaking the anti-adjacency between the pancakes $2k$ and $2k-1$.
    Consider now the clan $[n-1, n-2, \dots, 2]$ of size $n-2$. By 
		Lemma~\ref{clan of size 3}, it has to be split into clans of size at 
		most~$2$ before the second flip $(n)$, otherwise, it would not be possible to
		finish the sorting with only wastes. As only $\frac{n}{2}-2$ wastes remain to perform this splitting of the clans into clans of size $2$, the only solution
		is to split exactly between the pancakes $2k$ and $2k-1$ for each $2 \le k \le
		\frac{n}{2}-1$.
\end{proofclaim}

We can now prove that the method presented in this paper cannot provide a $\left(\frac{3}{2}n+1\right)$-flipping sequence that sorts $\overline{I_n}$ for even $n$.

\begin{theorem}
    Let $n$ be an even positive integer. If $T(n) = \frac{3}{2}n+1$, an optimal flipping 
		sequence of $-I_n$ cannot start with a flip $(n)$ (a waste) followed by all further 
		wastes, then a flip $(n)$ (an improve) followed by all further improves.
\end{theorem}

\begin{proof}
    Suppose by contradiction that such a flipping sequence exists. The proof is a case distinction according to the second flip.
    \begin{itemize}
				\item If the second flip is a flip $(n-1)$, the resulting pancake stack 
				is $[-2, -3, \dots, -n,1]$. By Claim~\ref{claim:breakclan2k-2k-1}, all 
				the following wastes cannot be flips $(n-1)$, otherwise they would not 
				break a clan. Hence, it is not possible to have the pancake $-n$ at the 
				top of the stack after all the wastes, contradicting 
				Claim~\ref{claim:secondn-flip-improve}.
				\item Otherwise, the second flip is a flip $(2k-1)$ for some $1\le k \le 
				\frac{n}{2}-1$, creating the substack $[-n,n-2k+1,n-2k]$. By 
				Claim~\ref{claim:breakclan2k-2k-1}, this substack cannot be broken by a 
				waste (but it can be fully flipped). Hence, after the last waste, the 
				stack has to start with $[-n,n-2k+1,n-2k]$, to ensure that the flip  
				$(n)$ is an improve, and that this substack has not been split. After the 
				flip $(n)$, the stack becomes $[-1, X, -(n-2k), -(n-2k+1),n]$ 
				for some 
				stack $X$, and has to be sorted with only improves. Furthermore, 
				by Claim~\ref{claim:breakclan2k-2k-1}, for all $ 1 \le i \le \frac{n}{2}-1$ the clan $(2i+1, 2i)$ 
				has not 
				been split. 
                
                Moreover, since all the remaining flips have to be improves, we can prove by induction that after the $j$-th improve after the second flip $(n)$, the first pancake of the stack is $-(2j-1)$ and all the clans $(2i+1, 2i)$ for $j \le i \le \frac{n}{2} -1$ have not been split. 
                
                Indeed, by hypothesis the first improve is the flip $(n)$ which puts the pancake $-1$ to the top of the stack and no clan $(2i+1, 2i)$ 
                for $ 1 \le i \le
                \frac{n}{2}-1$
                has been split. 
                
                Then suppose that after the $j$-th improve, the pancake at the
                top of the stack is $-(2j-1)$ and that all the clans $(2i+1,
                2i)$ for $j \le i \le \frac{n}{2} -1$ have not been split. The
                $(j+1)$-th improve has to set the pancake $-(2j-1)$ next to the
                pancake $2j$. Hence, since the clan $(2j+1,2j)$ has not been
                split, and since the pancake $2j$ has to be burnt face-down for
                this flip to be an improve, necessarily, this clan is in that
                order in the stack and thus the $j$-th improve is a flip between
                these two pancakes, putting the pancake $-(2j+1)$ to the top
                of
                the stack and the only clan than has been split is the clan
                $(2j+1,2j)$, which concludes the induction.
                
                By construction, for 
				$1 \le j \le \frac{n-2k}{2}$, the $j$-th improve after the second flip $(n)$ cannot modify the last three pancakes of the stack. 
				Hence, the $\left(\frac{n-2k}{2}\right)$-th improve after the second flip $(n)$ puts the pancake 
				$-(n-2k-1)$ to the top of the stack. This is a contradiction, as the next 
				flip has to be an improve and thus put it adjacent to the pancake 
				$n-2k$, which currently has the wrong side up in the stack. 
                \qedhere
    \end{itemize}
\end{proof}

\section{Conclusions and Future work}
\label{sec_conc}
Burnt pancake sorting is a classical problem going back 
to the work of~\cite{GP79}. The original
problem is to find lower and upper bounds for the
number of flips which are required to sort an arbitrary
stack of burnt pancakes of size $n$. As $ \overline{I_n}$
has been shown to belong to the worst-case instances
for all $ 1 \le n \le 22 $, except for $n=15$~\cite{CB95,Bou16}, and as it has
a very easy and natural structure, it is very
interesting to compute upper and lower bounds
for the number of flips which this instance needs 
to be sorted, called $T(n)$.

\cite{HS97}~gave an upper bound for $T(n)$
for $ n\equiv 3 \bmod 4$,
and~\cite{Cib11} a global lower bound, 
which matches the upper bound 
of~\cite{HS97} 
for $ n\equiv 3 \bmod 4 $.
\cite{Bou16}~computed the  values $T(2),T(3),\dots, T(27)$. 
\footnote{We have recomputed the values of $T(24)$ and 
$ T(26)$ to be $37$
and $40$, instead of being
$38$ and $41$, respectively, as mentioned in~\cite{Bou16}. 
On the other hand, we have confirmed the correctness of all
other values up to $27$, for $n$ odd
by our main result and 
for $n$ even by a BnB
program in C++.}
Apart from that, we are not aware of any other 
general result about $T(n)$ since then.
In this work, we have presented an upper bound 
for $ n\equiv 1 \bmod 4 $, which again 
matches the lower bound of~\cite{Cib11} .

Although the problem appears quite
simple at first glance,
the size of the burnt pancake graph and the complexity of finding optimal
strategies
make efficient computation highly challenging. Consequently, the tools developed in this work
to tackle this problem are of independent interest, beyond the (burnt) pancake sorting problem itself.

For future work we suggest the following:

\begin{itemize}
\item
The proofs of this work contains three
sub-cases, namely
$ n\equiv 1 \bmod 12 $,
$ n\equiv 5 \bmod 12 $ and
$ n\equiv 9 \bmod 12 $,
which have similar, long and technical proofs.
To shorten them it would be interesting to
find unified flipping sequences and corresponding
proofs for $ n \equiv 1 \bmod 4$ or even for all
odd $n$.

\item
The only open case is $T(n)$ for even $n$, where only two
values are possible,
namely $ \frac{3}{2}n+1$
and $ \frac{3}{2}n+2$.
In this context,
we have already shown that 
there are at least two further values
(except the case $n=2$) where
the lower bound is attained,
namely $n=24$ and $n=26$.

\item
Another problem would be to investigate whether $n=15$
is the only positive integer, where $\overline{I_n}$ is not a worst-case instance. Note that whenever $\overline{I_n}$ is a worst-case instance for some integer $n$, this work ensures that the diameter of the pancake graph is either $\frac{3}{2}n+1$ or $\frac{3}{2}n+2$. 
\end{itemize}

\section*{Acknowledgements}

This work was supported by the Kempe Foundation Grant No. JCSMK24-515 (Sweden).

This research was conducted using the resources of High Performance Computing Center North (HPC2N).The computations were enabled by resources provided by the National Academic Infrastructure for Supercomputing in Sweden (NAISS), partially funded by the Swedish Research Council through grant agreement no. 2022-06725.

\bibliography{aaai2026}

\end{document}